\documentclass[12pt]{amsart}
\newcommand {\supplus}{\mathop{{\supset}\llap{\raise 
0.5pt\hbox{\normalfont\small+}\hskip 0.5pt}}} 
\newcommand {\subplus}{\mathop{{\subset}\llap{\raise 
0.5pt\hbox{\normalfont\small+}\hskip 0.5pt}}}  
\newcommand {\Cee}    {{\mathbb  C}}

\newcommand {\Nee}    {{\mathbb  N}}

\newcommand {\Zee}    {{\mathbb  Z}}

\newcommand {\fg}     {{\mathfrak{g}}}    %
\newcommand {\fh}     {{\mathfrak{h}}}

\newcommand {\fo}     {{\mathfrak{o}}}
\newcommand {\fosp}   {{\mathfrak{osp}}}

\newcommand {\fS}     {{\mathfrak{S}}}

\newcommand {\fsl}    {{\mathfrak{sl}}}

\newcommand {\cal} {\mathcal}

\def \opname#1#2%
  {\expandafter\newcommand \csname #1\endcsname {{\mathop{#2}\nolimits}}}

\newcommand{\rmname}[1]
  {\expandafter\newcommand \csname #1\endcsname {{\operatorname{#1}}}}

\newcommand{\rmnameii}[2]
  {\expandafter\newcommand \csname #1\endcsname {{\operatorname{#2}}}}

\rmname{act}
\rmname{Ad}
\rmname{Add}
\rmname{ad}
\rmname{Alt}
\rmname{alt}
\rmname{Ann}
\rmname{antidiag}
\rmname{Ber}
\rmname{ber}
\rmname{Br}
\rmname{card}
\rmname{ch}
\rmname{Char}
\rmname{cem}
\rmname{cj}
\rmname{Cliff}
\rmname{cntr}
\rmname{codim}
\rmname{coind}
\rmname{const}
\rmname{col}
\rmname{cork}
\rmname{cpr}
\rmname{diag}
\rmnameii{Div}{div}
\rmname{Def}
\rmname{Der}
\rmname{Dim}
\rmname{End}
\rmname{Even}
\rmname{Ext}
\rmname{gr}
\rmname{Hom}
\rmname{HT}
\rmnameii{Ht}{ht}
\rmname{hwt}
\rmname{Id}
\rmname{id}
\rmname{ind}
\rmname{Ind}
\rmname{Inf}
\rmname{irr}
\rmname{Le}
\rmname{Lie}
\rmname{lwt}
\rmname{mult}
\rmname{Mor}
\rmname{nm}
\rmname{Ob}
\rmname{Odd}
\rmname{Osc}
\rmname{per}
\rmname{Pic}
\rmname{pr}
\rmname{pro}
\rmname{Prime}
\rmname{Proj}
\rmname{prt}
\rmname{pt}
\rmname{Q}
\rmname{qet}
\rmname{qtr}
\rmname{rd}
\rmname{rk}
\rmname{row}
\rmname{Res}
\rmname{salt}
\rmname{Sch}
\rmname{sch}
\rmname{SBr}
\rmname{scalar}
\rmname{Ser}
\rmname{sign}
\rmname{Smbl}
\rmname{spin}
\rmname{ssym}
\rmname{str}
\rmname{st}
\rmname{sgn}
\rmname{sq}
\rmname{symm}
\rmname{supp}
\rmname{Supp}
\rmname{St}
\rmname{Spec}
\rmname{Spm}
\rmname{tr}
\rmname{vpt}
\rmname{weyl}
\rmname{Weyl}
\rmname{Witt}

\opname{vvol}  {{v\hspace{-0.1ex}o\hspace{-0.02ex}l\/}}
\opname{pnt}  {\text{\normalfont pt}}
\opname{Span} {{Span}}
\opname{slim} {\overline{\lim}}
\opname{Vol}  {{V\hspace{-0.55ex}o\hspace{-0.02ex}l\/}}
\opname{Par}  {{P\hspace{-0.3ex}a\hspace{-0.05ex}r\/}}

\rmname{Mat}
\rmname{Bil}
\rmname{Diff}
\rmname{Ker}
\rmname{Herm}
\rmname{Coker}
\rmname{Conn}
\rmname{Covect}
\rmname{Vect}
\rmname{Int}

\rmnameii {IM} {Im}
\rmnameii {RE} {Re}

\opname{Aut} {{A\hspace{-0.2ex}u\hspace{-0.1ex}t\/}}
\opname{GL} {{G\hspace{-0.3ex}L}}
\opname{SL} {{S\hspace{-0.3ex}L}}
\opname{Exp} {{E\hspace{-0.2ex}x\hspace{-0.1ex}p\/}}
\opname{GQ} {{G\hspace{-0.2ex}Q}}
\opname{OSp} {{O\hspace{-0.25ex}S\hspace{-0.15ex}p\/}}
\opname{Out} {{O\hspace{-0.25ex}u\hspace{-0.15ex}t\/}}
\opname{Spp} {{S\hspace{-0.2ex}p\/}}
\opname{SpO} {{S\hspace{-0.2ex}p\hspace{-0.02ex}O\/}}
\opname{Pe} {{P\hspace{-0.25ex}e\/}}
\opname{SPe} {{S\hspace{-0.25ex}P\hspace{-0.25ex}e\/}}
\opname{Spin} {{S\hspace{-0.25ex}p\hspace{-0.05ex}i\hspace{-0.1ex}n\/}}
\opname{Iso} {{I\hspace{-0.25ex}s\hspace{-0.1ex}o\/}}
\opname{SSPe} {{S\hspace{-0.25ex}S\hspace{-0.15ex}P\hspace{-0.25ex}e\/}}
\opname{PeU} {{P\hspace{-0.25ex}e\hspace{-0.1ex}U\/}}
\opname{QU} {{Q\hspace{-0.15ex}U\/}}
\opname{U} {{U\/}}

\opname{cGQ} {{\cal G \hspace{-0.2em} Q \/}}
\opname{cSL} {{\cal S \hspace{-0.2em} L \/}}
\opname{cGL} {{\cal G \hspace{-0.2em} L \/}}
\opname{cOSp} {{\cal O \hspace{-0.2em} S \hspace{-0.3em} \it p\/}}
\opname{cPe} {{\cal P \hspace{-1.5pt} \it e\/}}
\opname{cVect} {{\cal V \hspace{-1.5pt} \it e\hspace{-0.1ex}c\hspace{-0.1ex}t\/
}}
\opname{cVol} {{\cal V \hspace{-1.5pt} \it o\hspace{-0.1ex}l\/}}
\opname{cAut} {{\cal A \hspace{-0.2em} \it u\hspace{-0.1em}t\/}}
\opname{cCovect} {{\cal C \hspace{-1.5pt}
     \it o\hspace{-0.1ex}v\hspace{-0.1ex}e\hspace{-0.1ex}c\hspace{-0.1ex}t\/}}
\opname{CW} {{C\hspace{-0.15ex}W}}

\newcommand {\ev} {{\bar0}}
\newcommand {\od} {{\bar1}}
\newcommand {\eps} {\varepsilon}

\newcommand {\tto} {\longrightarrow}

\newcommand {\bcdot}   {\mathbin{\hbox{\raise.4ex\hbox{\bf.}}}} % bold \cdot

\newcommand {\secno} {}
\newcommand {\ssecfont} {\normalfont\bf}

\newtheorem{Theorem}{\secno Theorem}

\newenvironment {th*}[1]
    {\gdef\thname{#1} \begin{thn}}%
    {\end{thn}}
\newtheorem{thn}[Theorem] {\thname}

\theoremstyle{definition}

\newenvironment {ex*}[1]
    {\gdef\thname{#1} \begin{exn}}%
    {\end{exn}}
\newtheorem{exn}[Theorem]{\thname}

\theoremstyle{remark}

\newtheorem{Remark}[Theorem]{\secno Remark}

\newenvironment {rem*}[1]
    {\gdef\thname{#1} \begin{remn}}%
    {\end{remn}}
\newtheorem{remn}[Theorem]{\thname}

\newcommand {\ssec}{\subsection*}

\newcommand {\ssbegin}[2]
  {\def \secno {\gdef \secno {}{\ssecfont #1. }}%
   \begin{#2}}
\begin{document}
\title{Orthogonal polynomials and Lie superalgebras}

\author{Alexander Sergeev} 

\address{Dept.  of Math., Univ.  of Stockholm, Roslagsv.  101, 
Kr\"aftriket hus 6, S-106 91, Stockholm, Sweden (On leave of absence 
from Balakovo Inst.  of Technology Technique and Control)\\ e-mail: 
mleites@matematik.su.se subject: for Sergeev}

\thanks{I am thankful to D. Leites, who raised the problem, for support 
and help.}

\begin{abstract} For $\fo(2n+1)$, in addition to the conventional set 
of orthogonal polynomials, another set is produced with the help of 
the Lie superalgebra $\fosp(1|2n)$.  Difficulties related with 
expression of Dyson's constant for the Lie superalgebras are 
discussed.
\end{abstract}



\maketitle

\section*{\S 0.  Introduction} 

\ssec{0.1.  History} In 1962 while studying statistical mechanics 
Dyson \cite{D} considered the constant term in the expression
$$
\mathop{\prod}\limits_{i\neq j}(1-\frac{x_{i}}{x_{j}})\text{ for }k\in \Nee,
$$
depending on indeterminates $x_{1}$, \dots , $x_{n}$.  Dyson 
conjectured the explicit form of this constant term.  His conjecture 
was soon related with the root system of $\fsl(n)$, generalized to 
other root systems of simple Lie algebras and proved.  The expressions 
obtained for the {\it Dyson constant} are called {\it Macdonald's 
identities}, see \cite{M}.

Let us briefly recall the main results. Let $\fg$ be a simple (finite 
dimensional) Lie algebra, $R$ its root system, $P$ the group of 
weights; $A=\Cee[P]$ the group of formal exponents of the form 
$e^\lambda$, where $\lambda\in P$; let $W$ be the Weyl group of $\fg$ 
and
$$
\triangle =\frac{1}{|W|}\mathop{\prod}\limits_{\alpha \in R}(1-e^\alpha).
$$
On $A$, define the scalar product by setting
$$
(f, g)=[f\bar g\triangle]_{0},\eqno{(0.0)}
$$
where $\overline{e^\lambda}=e^{-\lambda}$ and $[f]_{0}$ is the constant 
term of the power series $f$.

It turns out that

{\sl the characters $\chi_{\lambda}$ of finite dimensional irreducible
representations of $\fg$ are uniquely determined by their properties

{\em (i)} to form an orthogonal {\em (with respect to the form 
$(0.0)$)} basis in $A^{W}$, the algebra of $W$-invariant elements of 
$A$;

{\em (ii)} $\chi_{\lambda}=e^\lambda+$ terms with exponents $\mu$ such 
that $\mu<\lambda$}.

Here for Lie superalgebras I consider the following problem: what are 
the analogs of the scaral product $(0.0)$ (hence, of $\triangle$ and 
$W$) for which (i) and (ii) hold? If (i) and (ii) do not hold {\it as 
stated}, how to modify the definitions and the statement to make them 
reasonably interesting? 

\ssec{0.2. Main result} So far, there is not much that can count as 
a result, actually. I consider this note as a remark on the results 
from \cite{M} and a report on the work in progress. 

It turns out that for Lie superalgebras there is no function $\Delta$ 
(understood as a formal distribution) such that the characters of 
irreducible representations would satisfy (0.1), i.e., were orthonormal.
With one exception: the series $\fosp(1|2n)$.  Thanks to this 
excepton, the main results of this note are:

1) For $\fosp(1|2n)$ I reproduce an observation of Rittenberg and 
Scheunert \cite{RS} on a correspondence between irreducible 
$\fosp(1|2n)$-modules and $\fo(2n+1)$-modules.  (I also give a short 
and lucid demonstration of this correspondence.  \footnote{Leites 
informed me, that this demonstration basically coincides with the one 
Deligne communicated to Leites in 1991 (unpublished).}) From this 
correspondence I deduce in the $\fo(2n+1)$ case the existence of 
another set of orthogonal polynomials in addition to the set described 
in \cite{M}.

2) For any simple Lie superalgebra $\fg$ I can produce a function $\Delta$ 
for which the characters of the typical representations are 
orthonormal with respect to (0.0). I hope to return to this topic 
elsewhere.

\begin{Remark} Observe that for the simple Lie algebras, $\Delta$ 
appears in the Weyl integration formula: if $f$ is a class function on 
a compact Lie group $G$ such that $Lie(G)\otimes \Cee\cong \fg$ and 
$T\subset G$ is a maximal torus, then
$$
\int_{G}f~dg=\int_{T}f\Delta dt.\eqno{(0.2)}
$$
For the general Lie superalgebras the analog of identity (0.2) is 
unknown to me. Here are several little problems: not every simple Lie 
superalgebra (supergroup) over $\Cee$ has a compact form; the volume 
of those that have may vanish identically, cf. \cite{B}. 
\end{Remark} 

\section*[\S 1.  The orthogonality of $\fosp(1|2n)$-modules]{\S 1.  
The orthogonality of the characters of $\fosp(1|2n)$-modules} 

We recall some basic facts from the representation theory of 
$\fosp(1|2n)$ (see, e.g., \cite{K}) and (for convenience) $\fo(2n+1)$.

\ssec{1.1.  $\fosp(1|2n)$, its roots and characters} Set
$$
\renewcommand{\arraystretch}{1.4}
\begin{array}{l}
R_{\ev}=\{\pm\eps_{i}\pm\eps_{j}\text{ for }i\neq j; \; \pm 
2\eps_{i}\}, \; \; R_{\od}=\{\pm\eps_{i}\}; \\
S_{\ev}=\{\pm\eps_{i}\pm\eps_{j}\text{ for }i\neq j\}\subset R_{\ev};\\
2\rho_{0}=\sum_{i<j}(\eps_{i}-\eps_{j}+\eps_{i}+\eps_{j})+\sum_{i}2\eps_{i}=
2\sum_{i<j}\eps_{i}+2\sum_{i}\eps_{i}=2\sum(n-i+1)\eps_{i}.
\end{array}
$$

For the Lie superalgebras
$$
\rho=\rho_{0}-\rho_{1}, \text{ where } \rho_{1}=\frac12\sum_{i}\eps_{i}.
$$
All $\fosp(1|2n)$-modules are typical. The invariant bilinear form 
is $\str(\ad(x)^2)$. Explicitely, the restriction of this form onto 
Cartan subalgebra reads as follows:
$$
\renewcommand{\arraystretch}{1.4}
\begin{array}{l}
\str(\ad(x)^2)=\sum(\pm\eps_{i}\pm\eps_{j})^2+\sum(\pm 2\eps_{i})^2- 
\sum(\pm\eps_{i})^2=\\
2\sum(\eps_{i}\pm\eps_{j})^2+4\sum(\pm \eps_{i})^2-
\sum(\pm\eps_{i})^2=\\
2\sum\big((\eps_{i}+\eps_{j})^2+(\eps_{i}-\eps_{j})^2\big)+
6\sum(\eps_{i})^2=\\
(4n+2)\sum(\eps_{i})^2.
\end{array}
$$

The supercharacter of the finite dimensional irreducible 
$\fosp(1|2n)$-module $V^\lambda$ with highest weight $\lambda$ is 
$$
\sch V^\lambda= \frac{\mathop{\sum}\limits_{w\in W}\sign' (w)
e^{w(\lambda+\rho)}}{L},
$$
where
$$
L=\frac{\mathop{\prod}\limits_{\alpha\in 
R_{\ev}^+}(e^{\alpha/2}-e^{-\alpha/2})} 
{\mathop{\prod}\limits_{\beta\in R_{\od}^+}(e^{\beta/2}-e^{-\beta/2})}
$$
and $\sign'(w)$ is equal to $-1$ to the power equal to the number of 
reflections in the even roots $\alpha$ except those $\alpha$ for which 
$\frac{\alpha}{2}\in R_{\od}$.

The Weyl group $W$ of $\fosp(1|2n)$ is equal to 
$\fS_{n}\circ(\Zee/2)^{n}$ and
$$
\sign'(\sigma\cdot \tau_{1}\dots \tau_{n})=\sign(\sigma) \text{ for any
$\sigma\in \fS_{n}$ and $\tau_{i}$ from the $i$-th copy of $\Zee/2$}.
$$

Observe that 
$$
L=\sum_{w\in W}\sign'(w)e^{w\rho}.
$$
Indeed, apply the character formula to the trivial module.

In other words, everything is the same as for $\fo(2n+1)$ but instead 
of the character $\sign$ on $W$ we now take $\sign'$.

The unique, up to $W$-action system of simple roots in $\fg$ is of 
the form
$$
\Pi =\{\eps_{1}-\eps_{2}, \dots , \eps_{n-1}-\eps_{n}, \eps_{n}\}.
$$
Observe that
$$
\lambda+\rho=\sum(\lambda_{i}+n-i+\frac12)\eps_{i}.
$$
Since
$$
\renewcommand{\arraystretch}{1.2}
\begin{array}{l}
L=\mathop{\prod}\limits_{\alpha\in 
R_{\ev}^+}(e^{\alpha/2}-e^{-\alpha/2}) 
\displaystyle\frac{\mathop{\prod}\limits_{i}(e^{\eps_{i}}-e^{-\eps_{i}})} 
{\mathop{\prod}\limits_{i}(e^{\eps_{i}/2}-e^{-\eps_{i}/2})}=\\
\mathop{\prod}\limits_{\alpha\in S^+_{\ev}}(e^{\alpha/2}-e^{-\alpha/2})
\mathop{\prod}\limits_{\alpha\in R^+_{\od}}(e^{\alpha/2}+e^{-\alpha/2}).
\end{array}\eqno{(1.1)}
$$
We will use the latter expression of $L$ as well.

\ssec{1.2. $\fo(2n+1)$, its  roots and characters} Clearly, 
$R(\fo(2n+1))=S(\fosp(1|2n))_{\ev}\cup R(\fosp(1|2n))_{\od}$ and
$\rho$ is the half-sum of the positive roots; the restriction of the Killing 
form is proportional
to $\sum\eps_{i}^2$, the Weyl group is $W= \fS_{n}\circ(\Zee/2)^{n}$ 
and for the nontrivial homomorphism $\sign: \Zee/2\tto \{\pm 1\}$ we 
have
$$
\renewcommand{\arraystretch}{1.4}
\begin{array}{l}
\sign(\sigma\cdot \tau_{1}\dots \tau_{n})=\sign(\sigma) 
\sign(\tau_{1})\dots\sign(\tau_{n})\\
\text{ for any
$\sigma\in \fS_{n}$ and $\tau_{i}$ from the $i$-th copy of $\Zee/2$}.
\end{array}
$$
The system of simple roots is the same as for $\fosp(1|2n)$;
the character of the finite dimensional irreducible 
$\fo(2n+1)$-module $V^\lambda$ with highest weight $\lambda$ is
$$
\ch V^\lambda= \frac{\mathop{\sum}\limits_{w\in 
W}\sign(w)e^{w(\lambda+\rho)}}{L},
$$
where 
$$
L= \mathop{\prod}\limits_{\alpha\in 
S(\fosp(1|2n))^+_{\ev}}(e^{\alpha/2}-e^{-\alpha/2}) 
\mathop{\prod}\limits_{\alpha\in 
R(\fosp(1|2n))^+_{\od}}(e^{\alpha/2}-e^{-\alpha/2}).\eqno{(1.2)}
$$

Set
$$
(a, q)_{\infty}=(1-a)(1-aq)(1-aq^2)\dots
$$
Set, further (for an indeterminate $\epsilon$ such that $\epsilon^2= 
1$; it corresponds to the 1-dimensional odd superspace; we hope that 
the reader will not confuse $\epsilon$ with the root $\eps$):
$$
\triangle(q, t, \epsilon)=\frac{\mathop{\prod}\limits_{\alpha\in 
S(\fosp(1|2n))_{\ev}}(e^{\alpha}, q)_{\infty}}{(te^{\alpha}, 
q)_{\infty}}\cdot \frac{\mathop{\prod}\limits_{\alpha\in 
R(\fosp(1|2n))_{\od}}(\epsilon e^{\alpha}, q)_{\infty}}{(t\epsilon 
e^{\alpha}, q)_{\infty}}\; \; \text{ for $\alpha\in R(\fo(2n+1))$}.
$$

Let 
$$
t= q^k\text{ for }k\geq 0. 
$$
Then
$$
\renewcommand{\arraystretch}{1.2}
\begin{array}{l}
\triangle(q, t, \epsilon)=\mathop{\prod}\limits_{\alpha\in 
S(\fosp(1|2n))_{\ev}} \; \; 
\mathop{\prod}\limits_{r=0}^{k-1}(1-q^re^{\alpha})\cdot 
\mathop{\prod}\limits_{\alpha\in R(\fosp(1|2n))_{\od}} \; \; 
\mathop{\prod}\limits_{r=0}^{k-1}(1-\epsilon q^r e^{\alpha})=\cr 
\mathop{\prod}\limits_{\alpha\in R} \; \; 
\mathop{\prod}\limits_{r=0}^{k-1}(1-\epsilon^{p(\alpha)}q^re^{\alpha}),
\end{array}\eqno{(1.3)}
$$
where $p(\alpha)=0$ for $\alpha\in S(\fosp(1|2n))_{\ev}$ and $p(\alpha)=1$ for 
$\alpha\in R(\fosp(1|2n))_{\od}$.

Let $P$ be the group of weights of $\fosp(1|2)$ and let $A=\Cee[P]$ be 
the group of formal exponents of the form $e^\lambda$, where 
$\lambda\in P$.  Recall that $\lambda\in P$ if and only if 
$\lambda=\sum n_{i}\eps_{i}$, where $n_{i}\in\Zee$ for all the $i$.

The Weyl group $W=\fS_{n}\circ (\Zee/2)^n$ of $\fosp(1|2n)$ acts on 
$P$, hence, on $A$, as follows: $\fS_{n}$ permutes the $\eps_{i}$ and 
$(\Zee/2)^n$ changes their signs.

Observe that for $\fo(2n+1)$ the group of weights is larger than same 
for $\fosp(1|2n)$: the former includes the half-integer $n_{i}$.

\ssbegin{4.2}{Theorem} For $\fo(2n+1)$ there exists a unique (up to a 
constant factor) basis of $A^{W}$ consisting of $\lambda\in P^+$ such 
that

{\em a)} $F_{\lambda}=m_{\lambda}+ 
\mathop{\sum}\limits_{\mu<\lambda}u_{\lambda, \mu}m_{\mu}$, where
$m_{\lambda}=\mathop{\sum}\limits_{\nu\in\{\text{the orbit of } 
\lambda\}}e^{\nu}$ and the coefficients $u_{\lambda, \mu}$ are 
rational functions in $t$ and $\epsilon$;

{\em b)} for $f, g\in A$ define the pairing $(f, g)$ by means of 
formula $(0.0)$, where $\Delta$ is determined by $(1.3)$. Then
$$
(F_{\lambda}, F_{\mu})=0\text{ if } \lambda\neq\mu.
$$
\end{Theorem}

\begin{proof} Uniquness.  Since $F_{\lambda}=m_{\lambda}+ 
\mathop{\sum}\limits_{\mu<\lambda}u_{\lambda, \mu}m_{\mu}$, it follows 
that the transition matrix from $F_{\lambda}$ to $F_{\lambda}$ is an 
uppertriangular one, i.e., $m_{\lambda}=F_{\lambda}+ 
\mathop{\sum}\limits_{\mu<\lambda}a_{\lambda, \mu}F_{\mu}$.  If 
$F_{\lambda}{}'$ is another set of elements from $A$ with the 
properties needed, then the transition matrix from $m_{\lambda}$ to 
$F_{\lambda}{}'$ is also an uppertriangular one, so 
$F_{\lambda}{}'=\mathop{\sum}\limits_{\mu\leq\lambda}b_{\lambda, 
\mu}F_{\mu}$.  If both the bases are orthogonal, i.e., $(F_{\lambda}, 
F_{\mu})=(F_{\lambda}{}', F{}'_{\mu})=0$ for $\lambda\neq \mu$, this 
means that $F_{\lambda}=F_{\lambda}{}'\cdot C_{\lambda}$.

Existence. It suffices to prove the existence of an operator 
$D: A^{W}\tto A^{W}$ such that

i) $(Df, g)=(f, Dg)$;

ii) $Dm_{\lambda}=
\mathop{\sum}\limits_{\mu\leq \lambda}c_{\lambda, 
\mu}m_{\mu}$;

iii) if $\lambda, \mu\in P^+$ are distinct, then $c_{\lambda, \lambda} 
\neq c_{\mu, \mu}$.

Set
$$
\Delta_{+}=\mathop{\prod}\limits_{\alpha\in R^{+}}\; \; 
\mathop{\prod}\limits_{r=0}^{k-1}(1-\epsilon^{p(\alpha)}q^re^\alpha).
$$
In the Cartan subalgebra $\fh\subset\fg$ select $h_{1}=\eps_{1}^{*}$ 
which plays the role of a miniscule weight for the dual root system. 
More exactly, for $R^+=\{\eps_{i}-\eps_{j}; \eps_{i}+\eps_{j};  
\eps_{i}\}$ we have $\alpha (h_{1})=0$ or 1 for any $\alpha\in R^+$. 

Define the action of the operator $T=T_{h_{1}}$ by setting
$$
Te^\lambda=q^{\lambda(h_{1})}e^\lambda
$$
and define $D$ by setting
$$
Df=\sum_{w\in W}w(\Delta^{-1}_{+}T(\Delta_{+}f)).
$$
The operator $D$ is self-adjoint. Indeed,
$$
\triangle=w\triangle=w\triangle_{+}\overline{w\triangle_{+}}\text{ for 
any }w\in W.
$$

Further on,
$$
\renewcommand{\arraystretch}{1.6}
\begin{array}{l}
(Df, g) =\sum (w(\triangle_{+}^{-1}T(\triangle_{+}f)), g)=\\
\frac{1}{|W|}\mathop{\sum}\limits_{{w\in W}}\Big[w\big (
\frac{T(\triangle_{+}f)}{\triangle_{+}}\big)\cdot \bar g\cdot 
w\triangle_{+}\overline{w\triangle_{+}}\Big]_{0}=\\ 
\frac{1}{|W|}\mathop{\sum}\limits_{{w\in 
W}}\big[w(T(\triangle_{+}f))\cdot \bar g\cdot 
\overline{w\triangle_{+}}\big]_{0}= 
\frac{1}{|W|}\mathop{\sum}\limits_{{w\in 
W}}\big[w(T(\triangle_{+}f))\cdot \bar g w 
\overline{\triangle_{+}}\big]_{0}=\\
(\text{ since } g^{w}=g)=\frac{1}{|W|}\mathop{\sum}\limits_{{w\in 
W}}[w(T(\triangle_{+}f))\cdot \bar g 
\overline{\triangle_{+}}]_{0}=\\
(\text{ since the constant term is always $W$-invariant} )=
[T(\triangle_{+}f)\cdot \bar g \overline{\triangle_{+}}]_{0}.
\end{array}
$$
Similarly,
$$
(Dg, f) =
[T(\triangle_{+}g)\cdot \bar f \overline{\triangle_{+}}]_{0}.
$$
But, as is easy to verify,
$$
[T(f)\cdot \bar g ]_{0}=[T(g)\cdot \bar f ]_{0}.
$$
Therefore,
$$
(Df, g) =(f, Dg) \text{ for any }f, g\in A^W.
$$
Let us show that $D$ sends $A$ into $A$.  Set 
$\Phi=\Delta^{-1}_{+}T(\Delta_{+})$.  Then 
$D=\mathop{\sum}\limits_{w\in W}w(\Phi T(f))$ since $T$ is an 
automorphism.

Let us compute $\Phi$. Since $\alpha (h_{1})=0$ or 1 for any $\alpha\in 
R^+$, it follows that 
$$
T(e^\alpha)=\left\{\begin{matrix}qe^\alpha&\text{for }\;\alpha 
(h_{1})=1\\e^\alpha&\text{for }\;\alpha 
(h_{1})=0.\end{matrix}\right .
$$

The case $\epsilon=1$ is considered in [M1]. Therefore, in what follows 
we assume that $\epsilon=-1$. We have:
$$
\renewcommand{\arraystretch}{1.8}
\begin{array}{l}
\Phi=\frac{T(\Delta_{+})}{\Delta_{+}}= 
\displaystyle\frac{T(\mathop{\prod}\limits_{\alpha\in 
R^{+}}\; \; \mathop{\prod}\limits_{r=0}^{k-1}(1-\epsilon^{p(\alpha)} 
q^re^\alpha))} {\mathop{\prod}\limits_{\alpha\in 
R^{+}}\; \; \mathop{\prod}\limits_{r=0}^{k-1} 
(1-\epsilon^{p(\alpha)}q^re^\alpha)}=\\
\mathop{\prod}\limits_{\alpha\in 
R^{+}}\; \; \mathop{\prod}\limits_{r=0}^{k-1} 
\displaystyle\frac{1-\epsilon^{p(\alpha)}q^{\alpha(h_{1})}q^r 
e^\alpha}{1-\epsilon^{p(\alpha)}q^re^\alpha}=\\
\mathop{\prod}\limits_{\alpha\in R^{+}, \;
\alpha(h_{1})=1}\; \; \mathop{\prod}\limits_{r=0}^{k-1} 
\displaystyle\frac{1-\epsilon^{p(\alpha)}q^{r+1} 
e^\alpha}{1-\epsilon^{p(\alpha)}q^re^\alpha}=(t=q^r)\\
\mathop{\prod}\limits_{\alpha\in R^{+}, \; \alpha(h_{1})=1} 
\; \; \mathop{\prod}\limits_{r=0}^{k-1} 
\displaystyle\frac{1-\epsilon^{p(\alpha)}t 
e^\alpha}{1-\epsilon^{p(\alpha)}e^\alpha}=\\
\mathop{\prod}\limits_{\alpha\in 
R^{+}}\; \; \mathop{\prod}\limits_{r=0}^{k-1} 
\displaystyle\frac{1-\epsilon^{p(\alpha)}t^{\alpha(h_{1})} 
e^\alpha}{1-\epsilon^{p(\alpha)}e^\alpha}= 
\mathop{\prod}\limits_{\alpha\in 
R^{+}}\; \; \mathop{\prod}\limits_{r=0}^{k-1} t^{\alpha(h_{1})} 
\displaystyle\frac{\epsilon^{p(\alpha)}-t^{-\alpha(h_{1})} 
e^{-\alpha}}{\epsilon^{p(\alpha)}-e^{-\alpha}}=\\
\mathop{\prod}\limits_{\alpha\in 
R^{+}}\mathop{\prod}\limits_{r=0}^{k-1} t^{\alpha(h_{1})} 
\displaystyle\frac{1-\epsilon^{p(\alpha)}t^{-\alpha(h_{1})} 
e^{-\alpha}}{1-\epsilon^{p(\alpha)}e^{-\alpha}}=t^{\alpha(2\rho)} 
\mathop{\prod}\limits_{\alpha\in 
R^{+}}\mathop{\prod}\limits_{r=0}^{k-1} 
\displaystyle\frac{1-\epsilon^{p(\alpha)}t^{-\alpha(h_{1})} 
e^{-\alpha}}{1-\epsilon^{p(\alpha)}e^{-\alpha}}=\\
t^{\alpha(2\rho)} e^\rho\delta^{-1} \mathop{\prod}\limits_{\alpha\in 
R^{+}} (1-\epsilon^{p(\alpha)}t^{-\alpha(h_{1})} e^{-\alpha}),
\end{array}\eqno{(1.4)}
$$
where
$$
\delta=\sum_{w\in W}\sign'(w)e^{w\rho}.
$$
Observe that $w\delta= \sign'(w)$. Observe also that 
$$
\delta =\mathop{\prod}\limits_{\alpha\in 
R^{+}}(e^{\alpha/2}-\eps^{p(\alpha)} e^{-\alpha/2}).
$$
For any $X\subset R_{+}$ set
$$
\sigma(X)=\sum_{\alpha\in X}\alpha.
$$
If we simplify (1.4) by eliminating parentheses, then $\Phi$ takes the 
form
$$
\Phi=t^{2\rho(h_{1})} e^\rho\delta^{-1}
\sum_{X\subset R^{+}}\varphi_{X}(t)
e^{-\sigma(X)}),
$$
where
$$
\varphi_{X}(t)=(-1)^{|X|}\eps^{p(X)}t^{-\sigma(X)(h_{1})} 
$$
and $p(X)=\#(\text{odd roots that occure in $X$})$. Let us calculate 
$De^\mu$ for $\mu\in P$. Observe that
$w(Te^\mu)=q^{\mu(h_{1})}e^{w\mu}$,
hence,
$$
\renewcommand{\arraystretch}{1.6}
\begin{array}{l}
De^\mu=\delta^{-1}q^{\mu(h_{1})}t^{2\rho(h_{1})}
\mathop{\sum}\limits_{X\subset R^{+}}\varphi_{X}(t)
\mathop{\sum}\limits_{w\in W}\sign'(w)e^{w(\mu+\rho-\sigma(X))}=\\
q^{\mu(h_{1})}t^{2\rho(h_{1})} \mathop{\sum}\limits_{X\subset 
R^{+}}\varphi_{X}(t) \left( \delta^{-1}\mathop{\sum}\limits_{w\in 
W}\sign'(w)e^{w(\mu+\rho-\sigma(X))} \right )=\\
q^{\mu(h_{1})}t^{2\rho(h_{1})}
\mathop{\sum}\limits_{X\subset R^{+}}\varphi_{X}(t)
\chi_{\mu-\sigma(X)}.
\end{array}
$$
Here 
$$
\chi_{\mu-\sigma(X)}\neq 0\longleftrightarrow\text{ there exists a 
}\nu\in P^+\text{ such that }\nu+\rho=w(\mu-\sigma(X)+\rho).
$$
The last result does not depend on the choice of the one-dimensional 
character ($\sign$ or $\sign'$) on $W$ because the orbit of the weight 
$\lambda$ contains a dominant weight if and only if the 
$\lambda_{i}^2$ are pairwise distinct. Therefore,
$$
\begin{array}{l}
Dm_{\lambda}=\mathop{\sum}\limits_{\mu\in 
W_{\lambda}}De^\mu=t^{2\rho(h_{1})} \mathop{\sum}\limits_{X\subset 
R^{+}}\varphi_{X}(t) \mathop{\sum}\limits_{\mu\in W_{\lambda}} 
q^{\mu(h_{1})} \chi_{\mu-\sigma(X)}.
\end{array}\eqno{(1.5)}
$$
But
$$
\rho-\sigma(X)=\frac12\mathop{\sum}\limits_{\alpha\in
R^{+}}\eps_{\alpha}\alpha\eqno{(1.6)}
$$
where $\eps_{\alpha}=\pm 1$; hence, $w(\rho-\sigma(X))$ is of the same 
form (1.6) and, therefore,
$$
w(\rho-\sigma(X))=\rho-\sigma(Y)\text{ for some }Y\subset R^{+}.
$$
Thus, 
$$
\nu=w\mu-\sigma(Y)\leq w\mu\leq \lambda.
$$
Consequently, $Dm_{\lambda}= 
\mathop{\sum}\limits_{\nu\leq\lambda}C_{\lambda\nu}(h_{1})\cdot 
m_{\nu}$ and a) is verified.

Let us prove b).  In (1.4), $\nu=\lambda$ if and only if $Y=\emptyset$ 
and 
$$ 
w\mu=\lambda\Longleftrightarrow w(\rho-\sigma(X))=\rho,\; 
\mu=w^{-1}\lambda,\; X=R^{+}\cap(-wR^{+}).
$$
Therefore,
$$ 
C_{\lambda\lambda}(h_{1})=t^{2\rho(h_{1})} \sum_{X\subset R^{+}} 
q^{w^{-1}\lambda(h_{1})}\sign'(w) \chi_{R^{+}\cap(-wR^{+})}(t)
$$
because
$$ 
\renewcommand{\arraystretch}{1.4}
\begin{array}{l}
\chi_{R^{+}\cap(-wR^{+})}=(-1)^{|R^{+}\cap(-wR^{+})|} 
(-1)^{p(R^{+}\cap(-wR^{+}))}t^{(w^{-1}\rho-\rho)(h_{1})} =\\
\sign'(w) t^{(w^{-1}\rho-\rho)(h_{1})}.
\end{array}%\eqno{(1.4)}
$$
Hence, 
$$ 
C_{\lambda\lambda}(h_{1})=t^{\rho(h_{1})} \sum_{w\in W} 
q^{\lambda(wh_{1})}t^{\rho(wh_{1})},
$$
as in \cite{M}.  Therefore (see (6.14) in \cite{M}), if $\lambda\neq 
\mu$ and $\lambda,\mu\in P^+$, then $C_{\lambda\lambda}\neq 
C_{\mu\mu}$.
\end{proof}

\end{document}